\documentclass{article}[10pt]
\title{Product Gauges for Product Domains}
\author{Pat Muldowney}

\usepackage[pdftex]{graphicx}

\newtheorem{theorem}{Theorem}

\newtheorem{example}{Example}

\newcommand{\bL}{\mathbf{L}}
\newcommand{\bI}{\mathbf{I}}

\newcommand{\da}{\delta}
\newcommand{\vt}{\vspace{5pt}\\}
\newcommand{\R}{\mathbf{R}}

\newcommand{\bR}{\bar{\mathbf{R}}}

\newcommand{\D}{\mathcal{D}}

\newcommand{\E}{{\mathrm{{E}}}}

\newcommand{\N}{\mathcal{N}}
\newcommand{\I}{\mathbf{I}}

\newcommand{\bx}{\mathbf{x}}
\newcommand{\bX}{\mathbf{X}}

\newcommand{\sq}{\simeq}
\newcommand{\g}{\gamma}

\newcommand{\ve}{\varepsilon}

\newcommand{\bN}{\mathbf{N}}

\newcommand{\pr}{\mathbf{P}}

\newcommand{\proof}{\noindent\textbf{Proof. }\noindent}
\newcommand{\nproof}{\hfill$\mathbf{\bigcirc}\vspace{12pt}$ }

\date{}
\begin{document}
\maketitle
\begin{abstract}
In generalized Riemann integration (or Henstock-Kurzweil integration), the formation of the partitions used in Riemann sum construction is regulated by rules known as gauges. For instance, if $I$ is a one-dimensional real interval in a partition of the domain of integration, a term $f(x)|I|$ of a Riemann sum must satisfy a condition or rule $|I|<\delta(x)$ where $\delta(x)$ is a positive function (or gauge) defined for points $x$ of the domain. This article examines gauges for multi-dimensional domains $\Omega^T$ where $T$ can be infinite, and where $\Omega$ itself can be an infinite-dimensional Cartesian product.
\end{abstract}
\section{Introduction}
Product domains are a common theme throughout \cite{MTRV}, along with their gauges and integrals. The products involved are, essentially, $\R \times \R \times \R \times \cdots$, finitely or infinitely. The product element is $\R$, but other one-dimensional real intervals can easily be substituted.

In fact, any domain which lends itself to -complete integral construction can be used as a component of a product space, for formation of -complete integration in the latter. For instance, product spaces can themselves be component spaces of product spaces. So instead of $\R \times \R \times \R $, we could have $\R \times \R^3 \times \R^2$.

The latter, of course, is essentially the same as $\R \times \R \times \cdots \times \R, = \R^6$. What is the point of restructuring this product as $\R \times \R^3 \times \R^2$? This article addresses this and similar issues. 

The main point of the article is how to construct gauges for structured product domains.
Along the way, it is helpful to take a close look at the proof of Theorem 4 (divisibility theorem) in \cite{MTRV}.

\section{What is Three-dimensional Brownian Motion?}
In \cite{MTRV} a sample space $\Omega = \R^T$, usually with $T$ denoting an interval of time such as $]0, \tau]$, is used to represent stochastic processes $X_T$, as in
\[
X_T \sq x_T \left[\R^T, F_{X_T}\right].
\]
A sample path $x_T$ then consists of real numbers $x(t)$ (or $x_t$) for $0<t\leq \tau$,
\[
x_T = \left(x(t): 0<t\leq \tau\right)= \left(x_t: 0<t\leq \tau\right)
= \left(x_t\right)_{0<t\leq \tau}.
\]

For example, a Brownian motion process is represented as \[X_T \sq x_T \left[\R^T, G\right]\] where the distribution function $G(I[N])$ is defined on cells $I[N] \subseteq \R^T$. In this case the variable, or outcome, at time $t$ is $x(t)$, where the latter denotes the displacement in one dimension of the Brownian particle from the point of origin of the particle at time $t=0$.

So $x(0)=0$ is the origin or starting point of the particle, and $x(t)$ is the distance (at later time $t$) of the particle from the origin. This distance is unpredictable and indefinite, but is assumed to follow some rule of normal distribution; so it is ``observable'' in the sense used in \cite{MTRV}. In fact $x(t)$ is a possible outcome of a random variable $X_t$ (or $X(t)$). 

In other words, $x(t)$ is a ``sample value'' of the observable (or random variable) $X(t)$, and $x=(x(t))_{t\in T}$ is a ``sample path'' of the process (joint observable, joint random variable $X=(X(t))_{t\in T}$.
The collection $X_T$ ($=X$) of such $X_t$ (or $X(t)$), subject to a joint probability distribution 
\[
G(I[N]) = \pr \left[x(t) \in I(t), \;t \in N\right],
\]
gives the Brownian process $X_T \sq x_T \left[\R^T, F_{X_T}\right]$. 

As it stands, however, this mathematical formulation describes  the random particle in one direction only. Whereas, in reality, such a particle is moving randomly in three dimensions, not one.

\begin{example}\label{3-dim BM}
In order to have a mathematical representation of Brownian particle motion in three dimensions, the following scheme can be used. Suppose the motion is in mutually perpendicular spatial directions $s=1$, $s=2$, $s=3$ (that is, a horizontal dimension, a vertical dimension, and a third dimension of depth); and suppose the origin of the motion at time $t=0$ is $(0,0,0)$,
\[
\left(x(1,0), x(2,0), x(3,0)\right)
 =(0,0,0),\;\;\;\mbox{ or }\;\;\;\left(x_{1,0}, x_{2,0}, x_{3,0}\right) = (0,0,0).\]
Likewise, at time $t>0$, the particle displacement $x(s,t)$ in each of the three dimensions $s=1,2,3$ is
\begin{eqnarray*}
\left(x(1,t),x(2,t),x(3,t)\right) &=& \left(x(s,t):s=1,2,3;\;0<t\leq \tau\right) \vt
&=&
\left(x_{st}:s=1,2,3;\;0<t\leq \tau\right) \vt
\mbox{with }\;\;
X_{st}=X(s,t) &\sq &
x_{st} \left[ \Omega,G\right]
\end{eqnarray*}
 for $s=1,2,3;\;0<t\leq \tau$,
where sample space $\Omega = \R$ and distribution function $G$   satisfies (among other properties)
\[
G\left(I_{st}\right) = \pr \left[x_{st} \in I_{st}\right]
\]
for each cell $I_{st}$ in $\R$; the probability being given by the normal curve with mean zero and variance $t$, for $s=1,2,3$ and $0<t\leq \tau$. \nproof
\end{example}
Three-dimensional Brownian motion is not adequately described in Example \ref{3-dim BM} above. There are some further features which should be taken into account. (For instance, the actual probabilities governing the process have been left unspecified above.) At time $t$ a Brownian particle is located at position
\[
\bx(t)=\bx_t = \left( x(1,t),x(2,t),x(3,t)\right) 
=\left( x_{1t},x_{2t},x_{3t}\right) \in \R^3.
\]
For times $t_j$, $0=t_0<t_1<t_2< \cdots < t_{n-1}<t_n=\tau$, denote position  by
\[
\bx(t_j) =\bx_j = \left( x(1,t_j),x(2,t_j),x(3,t_j)\right) 
=\left( x_{1j},x_{2j},x_{3j}\right) =\left(x_{ij}\right)_{i=1,2,3}\in \R^3.
\]
Let $\bI =I_1\times I_2 \times I_3$ denote cells in $R^3$ where each $I_i$ ($i=1,2,3$) is a cell in $\R$.
For $1 \leq j \leq n$ let
\[
\bI_j = I_{1j}\times I_{2j} \times I_{3j}.
\]
For the Brownian particle moving in three dimensions the joint event $\bI_j$ ($1\leq j \leq n$) denotes the possibility that the particle position $\bx_j$ at times $t_j$ is in $\bI_j$; so
\[
x_{ij} \in I_{ij},\;\;\;i=1,2,3,\;\;\;j=1,2, \ldots n.
\]
The next step is to formulate a probability distribution, or probability values for the joint event $\bI_j$, $1 \leq j \leq n$; in other words, the probabilities
\[
\pr \left[\bx_j \in \bI_j\right]_{1\leq j \leq n},\;\;\;
=\;\;
\pr \left[x_{ij} \in I_{ij}\right]_{i=1,2,3,\;\;1\leq j \leq n}
\]
The conditions BM1 to BM7 (\cite{MTRV}, pages 305--306) are applicable to the increments $x_{ij}-x_{i,j-1}$ for $1 \leq j\leq n$ and for $i=1,2,3$, even though the increments are now expressed in all three dimensions of physical space rather than the one-dimensional simplification of \cite{MTRV}.

In particular, for any fixed $j$ (or $t_j$), the increments 
\[
x_{1j}-x_{1,j-1},\;\;\;\;
x_{2j}-x_{2,j-1},\;\;\;\;
x_{3j}-x_{3,j-1},
\]
in each of the three physical dimensions,
are each normally distributed random variables, statistically independent of each other and of the other increments at other times $t_{j'}\neq t_j$.

Thus, using the same reasoning as in the one-dimensional case in \cite{MTRV}, 
and writing $\bI = \bI_1\times \cdots \times \bI_n$, the probability distribution function $G$ for the three-dimensional Brownian motion is
\begin{equation}\label{3-dim G}
G(\bI)\;=\;\pr \left[\bx_j \in \bI_j\;:\;{1\leq j \leq n}\right]\;=\;
\prod_{i=1}^3
\left(
\prod_{j=1}^n
\left(
\frac{\int_{I_{ij}} \frac{-\left(x_{ij}-x_{i,j-1}\right)^2}{2\left(t_j-t_{j-1}\right)}dx_{ij}}
{\sqrt{2\pi\left(t_j-t_{j-1}\right)}}
\right)
\right). 
\end{equation}
This expression is a three-dimensional analogue of the one-dimensional construction in Section 6.8 (pages 284--288) of \cite{MTRV}. 

Remember, as discussed in pages 87--88 of \cite{MTRV}, a sample path $x_T$ can be thought of as a displacement-time graph \textbf{or} as a point of an infinite-dimensional Cartesian product space $\R^T$. But ``dimension'', in ``three-dimensional analogue'' above, refers, not to this issue, but to the difference between $\R$ (in $\R^T$) and $\R^3$ (in $\left(\R^3\right)^T$). In other words, it relates to the difference between $x_j -x_{j-1}$ and $x_{ij}-x_{i,j-1}$ ($i=1,2,3$).

The analysis of one-dimensional Brownian motion $X_T$ in \cite{MTRV} used the  stochastic process (defined in the -complete sense)
$
X_T \sq x_T[\R^T, G]
$
with distribution function $G$ defined in Section 6.8 of \cite{MTRV}. In the above discussion, some of the elements of a mathematical representation of three-dimensional Brownian motion have been introduced; in particular, a function $\bx_T$ which specifies the position in three-dimensional space, at each time $t$, of a Brownian particle,  and  a new version of a probability function $G(\bI)$ in (\ref{3-dim G}).

The question then arises,
whether it is  possible to build on these elements (that is, $\bx_T$ and $G(\bI)$) to form a -complete stochastic process
\[ 
 \bX_T \sq \bx_T \left[\Omega,G\right],
 \]
where the sample space $\Omega$ is the set $\left(\R^3\right)^T$ of all $\bx_T$.

To achieve this objective it is necessary to establish that a -complete integration system of gauge-constrained point-cell elements $(\bx, \bI)$ can be constructed in the domain $\Omega = (\R^3)^T$, in accordance with axioms DS1 to DS8 of  Section 4.1 in \cite{MTRV}; and in particular to establish existence of gauge-constrained divisions of domain $(\R^3)^T$, as in Theorem 4, pages 121--124 of \cite{MTRV}.

To prepare the ground for this task, some motivational examples follow.

\section{A Structured Cartesian Product Space}

\begin{example}\label{Bentropia}
Here is an example of multi-dimensional random variation which involves measurements which depend on parameters other than time $t$ and location in three-dimensional space. Suppose a single experiment on a moving object involves multiple measurements, with unpredictable outcomes, as follows:
\begin{enumerate}
\item
Spatial orientation $o$ of the object, given by two angles measured in radians; giving unpredictable outcomes $x_o(a_1,a_2)$ where
\[
0\leq x(a_1,o) = x_{1,o} <2\pi,\;\;\;\;0\leq x(a_2,o) = x_{2,o} <2\pi,
\]or $x_{i,o} \in S=[0,2\pi[,\;\;i=1,2$.
\item
Location $l$ of the object in 3-dimensional space, given by three distance values which measure (in centi\-metres, say) the co-ordinates of the location relative to some point of origin $(0,0,0)$; giving \[
x(d_i,l),\; = x_{i,l},\;\;\;({i=1,2,3})
\] where $x_{i,l} \in \R$ for $i=1,2,3$.
\item
Energy value  of the object (measured in joules, say); giving unpredictable outcome $x(e), =x_e$, where $x_e \geq 0$ (or $x_e \in \R_+$).
\end{enumerate} 
\end{example}
To summarize, write
\[
x_o=\left(x_{i,o}\right)_{i=1,2},\;\;\;x_l=\left(x_{i,l}\right)_{1=1,2,3},\;\;\;x_e=x_e,\;\;\;\bx=\left(x_o,x_l,x_e\right)
\]
where $ x_{i,o} \in S$, $x_{i,l} \in \R$ for each $i$, and $x_e \in \R_+$.
If appropriate distribution functions exist, this scenario can be formulated in  -complete terms, or Riemann-observable terms, (i.e.~as in \cite{MTRV}), as follows.

If outcomes $x_e$ (measurement of energy) have a probability distribution $F_{X_e}$ then $X_e \sq x_e[\R_+, F_{X_e}]$ is a basic observable.
If joint outcomes $x_l=x_{i,l}$ ($i=1,2,3$) have a joint probability distribution $F_{X_l}$ then
$X_l \sq x_l[\R^3, F_{X_l}]$ is a joint basic observable. 
If joint outcomes $x_o=x_{i,o}$ ($i=1,2$) have a joint probability distribution $F_{X_o}$ then
$X_o \sq x_o[S^2, F_{X_o}]$ is a joint basic observable.

However, the experiment consists of a single joint measurement of each of the variables. This suggests a single joint observable $X$ involving all of the parameters of the experiment. Let      
\begin{eqnarray}
&&M_o = \{1,2\},\;\;\;\;M_l=\{1,2,3\},\;\;\;\;M_e=\{1\},\;\;\;\;M=\{M_o,M_l,M_e\}\, ,\nonumber \vt
&&x_M = 
 \left((x_{1,o},x_{2,o}), (x_{1,l},x_{2,l},x_{3,l}), (x_e)\right) \;\in\; (S\times S) \times (\R\times \R \times \R) \times \R_+\;,\nonumber \vt
&& X_M \sq x_M\left[\Omega,\;F_{X_M}\right]\;\label{structure},
\end{eqnarray}
assuming a joint distribution function $F_{X_M}$ exists, defined on cells 
\[
I = I(M) = I(M_o) \times I(M_l) \times I(M_e)
\]
 of sample space 
\[
\Omega =(S\times S) \times (\R\times \R \times \R) \times \R_+.
\]
The symbols variously denoted by $M$ are labels corresponding to sets of dimensions $N$ in the theory formulated in \cite{MTRV} for domain $\R^T$.
The elements in (\ref{structure}) may appear somewhat over-structured. For instance, $x_M$ has three components consisting, respectively, of two, three and one joint measurements. 

It would  be  reasonable (and no doubt simpler) to treat the experiment as a whole as consisting of six joint measurements in unstructured sample space
$S\times S \times \R\times \R \times \R \times \R_+$ (or simply $\R^6$ with appropriate restrictions in the distribution functions); disregarding the 2,3,1 structure corresponding to measurement of angles, distances and energy, respectively.
However the more structured approach will be seen in due course to have some advantages.

In terms of the Riemann-observable (or -complete) theory in \cite{MTRV}, the \textit{joint basic} observable $x_M$ cannot be understood as a random variable; but a \textit{contingent} observable $f(X_M)$ may, under certain circumstances, qualify as a random variable. 

To illustrate, suppose the value of some imagined physical property, call it the \textit{bentropia} of the system, is measured as a real number  $\mathcal{B}$ obtained by means of a deterministic calculation $f$ on the orientation angles, the displacement distances, and the energy of the object,
\begin{equation}\label{bentropia}
\mathcal{B}=f(x_M) = f\left((x_{1,o},x_{2,o}), (x_{1,l},x_{2,l},x_{3,l}), (x_e)\right).
\end{equation}
This invention seeks to highlight a distinction between  joint-basic observables and  contingent observable:  \[x_M=(x_{o,1},x_{o,2}, x_{l,1},x_{l,2},x_{l,3}, x_e),\;\;\;\;f(x_M)=f(x_{o,1},x_{o,2}, x_{l,1},x_{l,2},x_{l,3}, x_e),\]
respectively. There are actual (as opposed to invented or meaningless) physical phenomena which lend themselves to this kind of mathematical description (\ref{bentropia}). But it is easier at this stage to avoid any complications which might arise with real (as opposed to invented) physical entities.

If $\int_\Omega f(x_M) F_{X_M}(I(M))$ exists then $\mathcal{B}$ is a contingent random variable, with expected value
\begin{eqnarray}
\E\left[\mathcal{B}\right]
&=&\int_\Omega f(x_M) F_{X_M}(I(M))
\label{structured integral}  \vt
&=&\int_{R_M}
f\left((x_{1,o},x_{2,o}), (x_{1,l},x_{2,l},x_{3,l}), (x_e)\right) F_{X_M}(I(M)), \nonumber 
\end{eqnarray}
where $R_M, =\Omega,$ denotes the structured product domain
\[(S\times S) \times (\R\times \R \times \R) \times \R_+\]
and $M=\{M_o,M_l,M_e\}$. 

For (\ref{structured integral}) to be meaningful in the Riemann-observable sense, $\int_{R_M}$ must qualify as Stieltjes-complete integral, as defined in Chapter 4 of \cite{MTRV}. Essentially, this means there must be gauges, divisions and Riemann sums which satisfy the conditions of Chapter 4. 

\textbf{The single most important step in this, from which most of the others follow directly, is the divisibility property which, in \cite{MTRV}, is established in Theorem 4, page 121.}

In this case,  (\ref{structured integral}) can simply be treated as an integral on $\R^6$---in fact a sub-domain of $\R^6$, with $x=(x_1,x_2,x_3,x_4,x_5,x_6)$ replacing $\bx=(x_o,x_l.x_e)$. But it is instructive to treat it instead as an integral on a Cartesian product of the three domains
\[
S\times S ,\;\;\;\;\;\; \R\times \R \times \R ,\;\;\mbox{ and }\;\; \R_+\;.
\]
The first two of these are themselves Cartesian products, and the associated elements in the three domains are, respectively,
\[
(x_o,I(M_o)),\;\;\;\;\;(x_l,I(M_l)),\;\;\;\;\;(x_e,I(M_e)),
\]
with gauges $\delta(x_o)$, $\delta(x_l)$, and $\delta(x_e)$.

The structured product domain  $\Omega =(S\times S) \times (\R\times \R \times \R) \times \R_+$ is to be denoted by $R_M$.
A point $\bx=x_M$ of $R_M$ is $(x_o,x_l,x_e)$. This is not quite the same as $(x_{1,o},x_{2,o},x_{1,l},x_{2,l},x_{3,l},x_{e}) \in \R^6$, though, as mentioned above, the latter is a bit simpler and works equally well in the theory.

Cells $I=I_M$ in $R_M$ are $I_o\times I_l \times I_e$. (Again, the more straightforward Cartesian product of six one-dimensional component cells can also be used in this case.)

The point-cell elements of $R_M$ are $(x_M,I_M)$, and these elements are associated if each of $x_o,x_l,x_e$ is associated with, respectively, $I_o,I_l,I_e$ in $S\times S$, $\R\times \R \times \R$, and $\R_+$. The point $x_M$ is an associated point (or tag-point) of cell $I_M$ if $x_M$ is a vertex of $I_M$. The set of tag-points of $R_M$ is $\bar R_M$, which is $R_M$ with ``points at infinity'' adjoined.

A gauge in $R_M$ is a function $\delta(x_M)$ defined for each tag-point of $R_M$. An associated point-cell pair $(x_M,I_M)$ is $\delta$-fine if each of the components
\[
(x_o,I_o),\;\;\;\;(x_l,I_l)\;\;\;\;(x_e,I_e)
\]
is $\delta$-fine in the manner described in \cite{MTRV} for each of the component domains $S\times S$, $\R \times \R \times \R$ and the domain $\R_+$. Of course, association and $\delta$-fineness can in this case be treated equally well by taking $R_M$ to be a straightforward six-dimensional Cartesian product; but the more structured approach may be more helpful in the longer run.

For this to be of any use, the conditions and properties of the integration theory of \cite{MTRV} Chapter 4 must be met by the structured elements $(x_M,I_M)$ of $R_M$. The single biggest issue here is Theorem 4 (pages 121--124 of \cite{MTRV}. We must prove $\delta$-fine divisibility of $R_M$ and its sub-domains. That is, given any gauge $\delta$, the integration domain $R_M$ has a $\delta$-fine division.

But, unlike Theorem 4, $R_M$ involves only finite Cartesian products. Therefore the proof of divisibility in $R_M$ is \textbf{not} like Theorem 4. Instead, it is similar to the standard or classical proof which is outlined on page 45 of \cite{MTRV}.
The well-known proof runs as follows.

\begin{theorem}
\label{divisibility of RM}
Given a gauge $\delta$ on $R_M$, there exists a $\delta$-fine division $\D$ of $R_M$, $\D = \{(x_M,I_M)\}$. Likewise if $R_M$ is replaced as domain by any cell or figure in $R_M$.
\end{theorem}

\proof
Assume, for contradiction, that the statement is false. Each component in the Cartesian product
\[
R_M= (S\times S)\times (\R\times \R \times \R)\times \R_+,
\]
can be bisected successively. At each stage a sub-domain is found to be non-divisible. As in Section 2.3 (page 45 of \cite{MTRV}), this gives a contradiction. \nproof

With this result to hand, it is not difficult to reproduce the theory of Riemann-observables for structured domains such as $R_M$. Therefore equation (\ref{structured integral}) is valid, and it is possible to incorporate this kind of random variable into the theory of \cite{MTRV}.

The preceding examples present some contrasting problems. In Example \ref{Bentropia} there are three ``dimensions'', labelled $o$, $l$, and $e$; representing, respectively, the two-, three-, and one-dimensional domains $S\times S$, $\R \times \R \times \R$, and $\R_+$.

In contrast, Example \ref{3-dim BM} has infinitely many dimensions, each labelled $t$ for $t \in ]0,\tau]$, and each dimension $t$ representing a three-dimensional Cartesian product space $\R^3=\R \times \R \times \R$, where, for $i=1,2,3$, the unpredictable measurement $x(i,t)$ (or $x_{i,t}$) is the $i$th co-ordinate of a Brownian particle at time $t$ ($0 <t \leq \tau$); so $x(t) =(x_{i,t})_{i=1,2,3} = (x_1(t), x_2(t), x_3(t) ) \in \R^3$.

Before proceeding with this construction, here are some background issues involving gauges for product spaces.

\section{Gauges for Product Spaces}
\begin{example}
\label{product gauge for Rn}
A gauge in a one-dimensional domain (such as $\R$) is a function $\da(y)>0$ defined for $y \in \bar \R$, $ = [-\infty, \infty]$. Section 2.17 (pages 79--81 of \cite{MTRV}) introduces gauges $\da$  for the finite-dimensional Cartesian product domain $\R^n$ where $n$ is any positive integer. Such a gauge $\da$ is simply $\da(x)>0$ for $x \in \bar \R^n$. Use the symbol $\da^a$ to distinguish such a gauge from $\da^b$ in (\ref{composite gauge for Rn}) below. 
To keep track of dimensions, denote the domain $\R^n$ by
\[
\R^n = \prod_{j=1}^n \R_j
\]
where $\R_j = \R$ for each $j$, and where a typical element of $\bar \R_j$ is denoted by $x_j$, and $x=(x_1, \ldots , x_n)$ is a typical element of $\bar \R^n$. Assume that a gauge $\da_j(x_j)>0$ is defined for each $\R_j$. The question posed here is whether a gauge 
\[
\da^b(x) = \da^b (x_1, \ldots , x_n)
\]
can be constructed from component gauges $\da_j$ for the product space $\R^n$  $= \prod_{1}^n \R_j$. In \cite{MTRV}, there is no presumption of elements $\da_j(x_j)$. 
Instead, a function $\da(x) = \da (x_1, \ldots , x_n)$ is defined for $x \in \bar \R^n$, without considering any of the one-dimens\-ion\-al elements $\da_j(x_j)$. However, it is perfectly feasible to construct a gauge $\da$ in $\R^n$ as the following composite of the $\da_j$:
\begin{equation}
\label{composite gauge for Rn}
\da^b (x) = \left( \da_1(x_1), \ldots , \da_n(x_n)\right),
\end{equation}
where component gauges $\da_j(x_j)$ are assumed to be given for each $\R_j = \R$. 

With this definition of gauge $\da^b$, the next question is to attribute meaning to ``$\da^b$-fine in $\R^n$''. Therefore suppose 
\[
(x,I), = \left((x_1, \ldots , x_n), I_1 \times \cdots \times I_n\right)
\]
is an associated pair in $\R^n$ (so $(x_j,I_j)$ are associated in $\R_j, = \R$, for $j=1, \ldots n$). Then, with gauge $\da^b$ given by (\ref{composite gauge for Rn}), we declare that $(x,I)$ is $\da^b$-fine in $\R^n$ if $(x_j,I_j)$ is $\da_j$-fine in $\R_j$ for $1 \leq j\leq n$.

Here are some properties of $\da^a$ and $\da^b$ gauges:

\begin{enumerate}
\item
Given $\da^b$, if $(x,I)$ is $\da^b$-fine then there exists a gauge $\da^a$ such that $(x,I)$ is $\da^a$-fine. For each $x = (x_1, \ldots , x_n$), simply 
take 
\[
\delta^a(x) = \min\left\{\delta_1(x_1), \ldots , \delta_n(x_n)\right\}.
\]
\item
Thus if $h(x,I)$ is defined on $\R^n$, and if $h$ is integrable on $\R^n$ with respect to $\da^b$ gauges (with integral value $\alpha$), then $h$ is integrable on $\R^n$ with respect to $\da^a$ gauges. For, by hypothesis, with $\ve>0$ given, there exists a gauge $\da^b$ so that, for every $\da^b$-fine division $\D^b$ of $\R^n$,
\[
\left|\alpha -(\D^b) \sum h(x,I) \right| < \ve.
\]
Now choose $\da^a$ as in 1, to complete the proof.
\item
Conversely, with $\da^a$ given, it is \textbf{not} possible to define a gauge $\da^b$ which ensures that every $\da^a$-fine pair $(x,I)$ is also $\da^b$-fine.
(This is analogous\footnote{The analogy: $\da^a(x)$ corresponds to $\da(y)$, and $\da^b(x)$ corresponds to constant $\da$. } to the extra ``discrimination'' that Riemann-complete gauge functions $\da(y)$ ($y \in \bar\R$) provide, in comparison with the constant $\da$ of standard Riemann integration.) To see this, consider any $\da^b$-gauge in a bounded two-dimensional domain $S=[0,1] \times [0,1]$,
\[
\da^b(x) = \da^b(x_1,x_2) = \left(\da_1(x_1), \da_2(x_2)\right)
\]
where $\da_j$  is a gauge in $[0,1]$ with $1>\da_j(x_j)>0$ ($j=1,2$). Let $\da^a(x)$ be
\[
\da^a(x) = \da^a(x_1,x_2) = \da_1(x_1)\da_2(x_2).
\]
Then, for each $x = (x_1,x_2) \in S$, there are cells $I_j$ ($j=1,2$) with
\[
\da^a(x) = \da_1(x_1)\da_2(x_2) < |I_j| <\da_j (x_j)
\]
for $j=1,2$, so $(x,I)$ is $\da^b$-fine, but not $\da^a$-fine.

\item
This means that---just like Riemann integration (with constant $\delta>0$) and Riemann-complete integration (with variable function $\delta(y)>0$)---there are functions $h(x,I), = h((x_1,x_2), I_1 \times I_2)$ which are integrable if the Riemann sums are formed with gauges $\da^a$, but are not integrable with gauges $\da^b$. 

We can say that $\da^a$-integration is \textbf{\emph{stronger}} than $\da^b$-integration, in the way that Riemann-complete integration is stronger than Riemann integration. (That is, a function $f$ which fails to be Riemann integrable may be Riemann-complete integrable.)

In review of this point, the role of a gauge (such as $\da^a$ or $\da^b$) is to restrict the associated pairs $(x,I)$ which can be admitted as members of division $\D$ in the formation of Riemann sums $(\D)\sum h(x,I)$ which are used to estimate $\alpha =\int h(x,I)$, with
\[
\left| \alpha - (\D) \sum h(x,I)\right| < \ve\]
as integrability criterion. 

\emph{\textbf{The more restrictive the gauge, the strong\-er the integral.}} 
\item
Taking this to extremes, if the gauge is such that \textbf{every} $(x,I)$ is admissible, then only\footnote{Theorem 67 (page 180 of \cite{MTRV}) gives an example of an integrand $h(x,I)$, taking the form $ f(x) \mathbf D(I)$, which is integrable if and only if $f$ is constant. } the constant function $h=\alpha$ is integrable. At the other extreme, if \textbf{no} $(x,I)$ is admissible, then \textbf{every} function $h$ is integrable. The mathematically useful integration scenarios lie between these two extremes.

\item
Despite the preceding issues, Fubini's theorem (\cite{MTRV}, pages 160--165) sometimes enables us to get a connection between $\da^a$-integrals and $\da^b$-integrals in $\R^n$. For instance, suppose $(x,I) = \left((x_1,x_2),I_1 \times I_2\right)$ in $\R^2$, and suppose
\[
h(x,I) = h_1(x_1,I_1)h_2(x_2,I_2).
\]
Now suppose $h$ is $\da^a$-integrable in $\R^2$. Then Fubini's theorem implies
\begin{eqnarray*}
\int_{\R^2} h(x,I) &= &\int_{\R_2} \left(\int_\R h_1(x_1,I_1) \right) h_2(x_2,I_2)\vt
& = &\int_\R h_1(x_1,I_1)\int_\R h_2(x_2,I_2).
\end{eqnarray*}
Since the latter two integrals are one-dimensional, the $\da^a$ gauges reduce to $\da^b$ gauges in these two integrals.
\item
The preceding argument based on Fubini's theorem (for $\da^a$-integrals) can be extended to integrals in $\R^n$; and it can be extended to some other kinds of integrands provided the iterated integrals of Fubini's theorem are expressible as successive integrals with respect to variables $x_j$ ($j=1, \ldots , n$), each of which is reducible to a $\da^b$-integral.

\item
\emph{\textbf{The $\da^a$-integral is stronger than the $\da^b$-integral. But a survey of the results and properties of the -complete integral $\int_{\R^n}h(x,I)$ ($\da^a$ version) in \cite{MTRV} shows that all results are equally valid for the corresponding statements expressed in terms of gauges $\da^b$.}}
\end{enumerate}

\end{example}
The definition of the gauge integral  is as follows. A function $h(x,I)$ is integrable on $\R^n$, with integral value $\alpha$ if, given $\ve>0$, there exists a gauge $\da$ such that, if $\D$ is a $\da$-fine division of $\R^n$, then
\[
\left| \alpha - (\D) \sum h(x,I)\right| < \ve.
\]
If there are no $\da$-fine divisions $\D$ of $\R^n$, then this definition is vacuous, implying (as per point 5.~above) that \textbf{every} function $h$ is integrable. Existence of $\da$-fine divisions is therefore the foundation of every form of -complete integration.

For gauges of form $\da^a$, this issue is addressed in various places in \cite{MTRV} (such as Section 2.3, page 45) using the method of successive bisection. Existence of divisions relative to gauges of form $\da^b$ can be easily established by the same method.
\begin{theorem}
\label{divisibility theorem for delta b in Rn}
Given a gauge $\da^b$ in $\R^n$, there exists a $\da^b$-fine division of $\R^n$.
\end{theorem}
\proof
Assuming non-divisibility (for contradiction), successive bisection delivers a point $x=(x_1, \ldots , x_n) \in \bar\R^n$, and, using the gauge values
$\da_j(x_j)$ for this point, a $\da^b$-fine cell $I$ can be found, and this proves divisibility of $\R^n$. \nproof

\section{Gauges for Infinite-dimensional Spaces}
The preceding discussion demonstrates that there are alternative ways of setting up a -complete system of integration for product spaces. These ideas are also helpful in formulating such a system when the product is infinite; that is, for domain ``$\R^\infty$'' rather than $\R^n$. This is illustrated as follows.

\begin{example}
\label{alternative gauge for RT}
Suppose $T=\,]0,\tau]$, and the domain of integration is the infinite-dimensional Cartesian product space
\[
\R^T = \prod_{t \in T} \R = \prod_{t \in T} \R_t
\]
where $\R_t = \R$ for each $t \in T$. In \cite{MTRV}, several approaches to this are discussed; but the method which is selected for use throughout the rest of the book is described in Section 4.2, pages 116--119. This can be summarized as follows.
\begin{equation}
\begin{array}{rrllrllll}
L: & \bar \R^T & \rightarrow & \N(T),&\;\;\;\;\;\; x_T &\mapsto & L(x_T) & \in& \N(T), 
\vt
\delta: & \bar \R^T \times \N(T) & \rightarrow& ]0, \infty[, &\;\;\;\;\;\;
(x_T, N)& \mapsto &\delta (x_T, N) &>&0.
\label{definition of gauge g 0}
\end{array}
\end{equation}
The purpose of $\delta$ is, with finite dimension-set $N\supset L(x_T)$ selected from the preceding line, to regulate the lengths $|I(t)|$ of the restricted edges $I(t)$ ($t \in N$) of associated pairs $(x_T, I[N])$ in $\R^T$.

Only a finite number of edges of $I[N]$ is restricted. As discussed in \cite{MTRV}, this is suggestive of gauge-restriction of edges of finite-dimensional cells $I(N)$ in the finite-dimensional Cartesian product domain $\R^N$. Example \ref{product gauge for Rn} describes two ways of providing gauge-restriction for $I(N)$:
\begin{enumerate}
\item[(a)]
A gauge $\da^a(x(N))>0$, with $|I(t)| < \da^a(x(N))$ for each $t \in N$;
\item[(b)]
A gauge $\da^b(x(N)) =\left(\da_t(x(t)): t\in N\right)$, with 
 $|I(t)| < \da^b(x(t))$ for each $t \in N$, where $\da^b(y)>0$ is a gauge in $\R=\R_t$ for each $t\in N$.
 \end{enumerate}
But while \cite{MTRV} mentions option (a) above, the book chooses, instead, the following form of gauge-restriction for $I[N]$:
\begin{enumerate}
\item[(c)]
A gauge $\da^c(x_T,N)>0$, with $|I(t)| < \da^c(x_T,N)$ for each $t \in N$.
\end{enumerate}
In other words, for each $t \in N$, the restriction on the edge $I(t)$ of $I[N]$ depends on \textbf{every} component $x(t)$ ($t \in T$, infinite) of $x(T) \in \bar\R^T$; and not just on the finite number of components $x(t)$ for $t \in N$, as described in (a) above.

It is easily seen that $\da^c$ gives a ``stronger'' integral than $\da^a$ (which is, in turn, stronger than $\da^b$, as discussed in Example \ref{product gauge for Rn}.) But the ``extra strength'' of $\da^c$ is superfluous for present purposes. Examination of the results in \cite{MTRV} (which use $\da^c$) show that they can also be obtained with gauges $\da^a$ or $\da^b$.

\emph{\textbf{For present purposes, it is more convenient here to use gauges $\da^b$ in $\R^T$, and also in other infinite-dimensional Cartesian product spaces under discussion.
}}  

Accordingly, and for purposes of reference, here is the relevant description of the chosen gauge for $\R^T$:

\begin{equation}\label{definition of gauge g 0 (b)}
\begin{array}{rrllrllll}
L: & \bar \R^T & \rightarrow & \N(T),&\;\;\;\;\;\; x_T &\mapsto & L(x_T) & \in& \N(T), \label{N for RT (b)}\vt
\delta_N: & \bar \R^N  & \rightarrow& \R_+^N, &\;\;\;\;\;\;
x (N)& \mapsto &\left(\delta_t(x(t))\right)_{t \in N} &\in & \R_+^N,
\label{delta for RT (b)}
\end{array}
\end{equation}
where $\R_+ =\,]0, \infty[$.
The gauge for $\R^T$ is the pair $(L, \da_N)$, which can be denoted by $\g^b$ since the $\da_N$ component has the form $\da^b$ of Example \ref{product gauge for Rn}. But whenever it is not likely to be confused with the gauges in \cite{MTRV}, the notation $\g=(L, \da_N)$ can be used.

The essential components of a gauge in $\R^T$ are therefore as follows.
\begin{itemize}
\item
For each $t \in T$, there is a gauge $\da_t$ in the corresponding domain $\R$.
\item
If $T$ contains only a \textbf{finite} number of elements $t$, then a gauge $\da$ in $\R^T$ is
\[
\left(\delta_t(x(t))\right)_{t \in T}
\]
where $(x_T, I_T)$ is $\da$-fine if $(x_t,I_t)$ is $\da_t$-fine for each $t \in T$.
\item
If $T$ contains \textbf{infinitely many} elements $t$, then, as in (\ref{definition of gauge g 0 (b)}) above, a  function $L$ is used to select, for each $x_T \in \bar\R^T$,
\[
L: \;\; \bar \R^T \; \rightarrow \; \N(T),\;\;\;\;\;\; x_T \;\mapsto \; L(x_T) \;\in \; \N(T).
\]
\item
For infinite $T$ a gauge $\g$ in $\R^T$ can then be written
\[
\g =\left( L, \left(\delta_t(x(t))\right)_{t \in T} \right),
\]
and an associated pair $(x,I[N])$ in $\R^T$ is $\g$-fine if
\[
N \supseteq L(x)\;\;\;\mbox{ and }\;\;\;(x_t,I_t)\;\mbox{is $\da_t$-fine for each}\;t \in N.
\]
\item
Thus, for infinite $T$, a gauge $\g$ can be written alternatively as
\[
\left(L, \da_N\right)\;\;\;\mbox{ or }\;\;\;\left(L, \left\{\da_t\right\}_{t \in T}\right).
\]
\end{itemize}

\end{example}
As divisibility  is the basis of -complete integration, $\g$-fine divisibility must be established as in Theorem \ref{divisibility of RM}. But since $T$ is infinite, the simple bisection argument of Section 2.3, page 45 of \cite{MTRV}, is not sufficient and a proof on the lines of Theorem 4 (pages 121--124 of \cite{MTRV}) is needed.

The formulation of the latter proof is, by necessity, burdened with quite a lot of technical notation. The following example seeks to give the underlying idea in a simpler, less technical way.

\begin{example}\label{divisibility finite dim}
Suppose $S=[0,1]$ and the domain of integration is 
\[
S^3=\,]0,1] \times \,]0,1] \times \,]0,1],
\]
 with associated elements \[
(x,I) =\left( (x_1,x_2,x_3), I_1 \times I_2 \times I_3\right)
\]
where, for $i=1,2,3$, $x_i$ is an element of $[0,1]$ and  $I_i$ is a cell $]u,v] \subseteq S=\,]0,1]$.
A gauge is a function $\delta (x_1, x_2, x_3)>0$ defined for $0\leq x_i \leq 1$, $(i=1,2,3)$. The usual way to prove $\delta$-divisibility is by a process of successive bisection:
\[
 \left ]\frac{q_{i_1}}{2^r}, \frac{q_{i_1}+1}{2^r}\right]
 \times 
\left ]\frac{q_{i_2}}{2^r}, \frac{q_{i_2}+1}{2^r}\right]
\times \left ]\frac{q_{i_3}}{2^r}, \frac{q_{i_3}+1}{2^r}\right],
\]
 $q_{i_k}= 0, 1, \ldots , 2^r-1$, $k=1,2,3$, and $r=1,2,3, \ldots$. Then an initial assumption of non-divisibility eventually produces a contradiction. But, in line with the proof of Theorem 4 (pages 121--124 of \cite{MTRV}), and dropping much of the technical notation, denote the three dimensions by $t_1,t_2,t_3$ instead of $1,2,3$. Think of the domain $S^3$ as a cubical block of cheese, each of whose edges has length 1. 
\begin{enumerate}
\item 
 Instead of bisecting the unit cube into successively smaller cubes, select one dimension only, $t_1$, for successive bisection. The initially cubical block of cheese is converted, by successively slicing along one dimension, into successively thinner slices of cheese, of which, at each stage, at least one slice, $J^r=I_{t_1}^r \times S \times S$, must be non-$\delta$-divisible from the initial assumption of non-divisibility.
As the bisection value $r$ tends to infinity a value $y_1$ (or $y_{t_1}$) is arrived at by the usual bisection argument, and, at each stage of bisection $y_{t_1}$ is the fixed $t_1$ co-ordinate of a point $(y_{t_1}, x_{t_2}, x_{t_3})$ which is the associated point or tag-point of slice $J^r$. 
\item
The next step is to consider a domain $S\times S$ with dimension $t_1$ removed. We take a paper-thin ``slice'' of cheese (actually so thin that it has thickness zero) consisting, in space, of the points $(y_{t_1}, x_{t_2}, x_{t_3})$ with $y_{t_1}$ fixed and $0 \leq x_{t_i} \leq 1$ for $i=1,2$. This is no longer really a slice of cheese; it is mathematically a two-dimensional Cartesian product domain $[0,1] \times [0,1]$. Define a gauge on this domain by
\[
\delta(x_2,x_3)=\delta(x_{t_2},x_{t_3}) := \delta(y_{t_1},x_{t_2},x_{t_3}).
\]
By the original assumption, this two-dimensional domain must be non-$delta$-divisible (with the new, two-dimensional meaning of $\delta$). Because if there is a $\delta$-fine division $\D$ of $S \times S$ (or, if preferred, of $\{y_{t_1}\} \times S \times S$, with
\begin{eqnarray*}
\D& =& \left\{\left(x_{t_2}, x_{t_3},\; I_{t_1} \times I_{t_2} \right) \right\}; \vt &&\mbox{or, alternatively, }\;  \left\{\left((y_{t_1},x_{t_2}, x_{t_3}),\; \{y_{t_1}\} \times I_{t_1} \times I_{t_2} \right) \right\},
\end{eqnarray*}
then $r$ can be chosen so that $2^{-r}$ is less
than $\delta\left(y_{t_1},x_{t_2}, x_{t_3}\right)$ for each of the finite number of tag-points $\left(y_{t_1},x_{t_2}, x_{t_3}\right)$ of $\D$ with fixed $y_{t_1}$. But then a $delta$-fine division of $I^r \times S \times S$ can be produced from the terms of $\D$:
 \[
\left\{\left((y_{t_1},x_{t_2}, x_{t_3}),\; I_{t_1}^r \times I_{t_1} \times I_{t_2} \right) \right\}.
\]
So if the two-dimensional ``paper-thin cheese sheet'' is divisible then for sufficiently large $r$ the ``somewhat thin'' (three-dimensional) cheese slices of 1 above are divisible, which  contradicts the original hypothesis of non-divisibility.
\item
Now select a second dimension $t_2$, and again commence successive bisection of the (vanishingly thin) ``cheese-sheet''; this time in the direction $t_2$, leaving direction $t_3$ unbisected. This produces a succession of non-divisible two-dimensional strips of width $2^{-r}$ and length 1, leading as before to a fixed value $y_{t_2}$ from the intersection of the strips.
\item
Finally, we get a non-divisible (one-dimensional) line in direction $t_3$, and further bisection gives fixed $y_{t_3}$, giving a fixed point $y=(y_{t_1}, y_{t_2}, y_{t_3})$ in $S\times S \times S$. Taking $r$ large enough, a cell $I_{t_1}^r \times I_{t_2}^r \times I_{t_3}^r $ is found to be divisible, with $y$ as tag-point, and this is found to contradict the original assumption of non-divisibility of $S\times S \times S$. 
\end{enumerate}
\end{example}
As pointed out at the start of Example \ref{divisibility finite dim} the exact same result can be obtained directly by simultaneously bisecting in all three dimensions.

So why take the rather more complicated route of bisecting in one dimension at a time? The next example demonstrates that this method can indeed be useful.

\begin{example}\label{divisibility infinite dim}
Instead of $S\times S \times S$ (with $S=[0,1]$), suppose $T$ represents a countably infinite number of dimensions $\{t_1,t_2, t_3, \ldots \}$, and suppose $R=S^T$ is the domain to be considered. In a finite-dimensional domain such as $S^3$ or $S^n$, a gauge is a function $\delta(x)>0$. But in infinite dimensions $\prod_t \left\{S: t \in T\right\}$, there is an additional condition. For each finite subset $N \subset T$, $\delta (x(N))>0$ is defined\footnote{Strictly speaking, \cite{MTRV} uses the more sophisticated function $\delta(x_T, N)>0$. But as mentioned in Example \ref{alternative gauge for RT}, $\delta(x(N))>0$ works fine. Also the version $\delta_N(x(N))=\left(\da_t(x_t):t \in N\right)$.} on the finite-dimensional domain $\prod_t \left\{S: t \in N\right\}$. In a way, this makes a supposedly infinite-dimensional $\int_{S^T}$ look like a straightforward finite-dimensional $\int_{S^N}$. 

But there is a catch. It is true that something like $\int_{S^N}$ occurs---in the sense of something resembling a Riemann sum estimate of $\int_{S^N}$. 
But there is a further condition on the gauges: for each $x_T$, a minimal finite set $L(x_T) \subset T$ is specified, and, in forming Riemann sum estimates of $\int_{S^T}$, we allow only those finite-dimensional cells $I(N)$ for which $N \supseteq L(x_T)$. The definition of a gauge $\g, = (L, \delta)$, is such that, in choosing gauges, the cardinality of the finite sets $N \supseteq L(x)$ can be arbitrarily large, while the corresponding numbers $\delta(x(N))$ can be arbitrarily small. These ideas are illustrated graphically in \cite{MTRV} pages 81, 87, and (especially) 102.

The issue  here is, with any gauge $\g=(L, \da)$ given, that the domain $S^T$ must be $\g$-divisible for the theory to work. There is a proof of this kind in \cite{MTRV} Theorem 4, pages 121--124. In Example \ref{divisibility finite dim} above, aspects of this proof are illustrated, but only for $T$ of finite cardinality. The underlying intent of the proof in Example \ref{divisibility finite dim} (in which $T$ is a finite set) is demonstrated in the present example in which $T$ is infinite.
\begin{enumerate}
\item
Assume $S^T$ is non-$\g$-divisible.
\item
As in Example \ref{divisibility finite dim} above, choose any dimension label. That is, choose any element of $T$. We have already enumerated $T$ as $\{t_1, t_2, t_3, \ldots \}$, so let the choice be $t_1$. 
\item 
There is no harm in continuing to visualise $S^T$ as a cube or ``block of cheese''. (Anyway, our powers of geometric visualisation do not easily extend
to more than three dimensions.) In reality $S^T$ is a hyper-cube, so to speak. That is, a hyper-block of cheese.
\item
As in Example \ref{divisibility finite dim}, bisect successively $S^T$ in dimension $t_1$ only. (That is, successively re-slice the hyper-block of cheese into ever-thinner slivers.) As before this produces cells (slices) $J^1 \supset J^2 \supset J^3 \supset \cdots$, each of which must be non-divisible by assumption 1 above. As before, a common, fixed number $y_{t_1}$ is arrived at.
\item
Now get rid of dimension $t_1$, so the hyper-block of cheese becomes a ``flat hyper-block'' or hyper-plane $S^{T \setminus \{t_1\}}$. Choose a second dimension $t_2$, bisect again, and find another fixed co-ordinate $y_{t_2}$.
\item
Thus far, the argument is similar to Example \ref{divisibility finite dim}. In that case we arrived at dimension $t_3$, bisected in that dimension, found $y_{t_3}$, and then stopped there because there were no more dimensions to slice. The proof was then completed by finding a contradiction to assumption 1. 
But in this case there are further dimensions $\{t_4, t_5, \ldots\}$, and we can \textbf{never} stop! In other words, for $j=1,2,3, \ldots$, we find (inductively) fixed co-ordinates $y_{t_j}$, giving (by induction) a single fixed point $y_T$ in $S^T$.
\item
Thus the final (or $t_3$) step in Example \ref{divisibility finite dim} is not available to us when $T$ is infinite. To see how to salvage the argument, remember that the gauge $\g$ is not just $\da$, but also includes a factor $L(x_T)$, a finite subset of $T$; a possibly different finite set for each $x_T$. Essentially, it is the finiteness of the set $L(y_T)$ (see step 6) that gives a contradiction to 1, forcing an eventual stop to the iteration. Return to steps 2, 3, and 4 above, where the fixed co-ordinate $y_{t_1}$ has been found. Now consider the hyperplane
\[
Z^1_T = \left\{z_T\right\} =\left\{(y_{t_1}, z_{t_2}, z_{t_3}, \ldots): z_{t_j} \in S,\; j=2,3,4, \ldots\right\}.
\]
\item
Suppose (for contradiction) that there is, in $Z^1_T$, a point \[z^1=(y_{t_1}, z_{t_2}, z_{t_3}, \ldots)\] for which 
\[
L(z^1) \cap \{t_2, t_3, \ldots \} = \emptyset,\;\;\mbox{ so }\;L(z^1) = \{t_1\}.
\]
Choose $r$ so that 

\begin{eqnarray}
2^{-r}& <& \da \left(z^1,\;\{t_1\}\right) \nonumber\vt
 \mbox{[ or } 
2^{-r} &< &\da \left(z^1(\{t_1\})\right),\nonumber\vt
 \mbox{or }
2^{-r}& < &\da \left(z^1_{t_1}(\{t_1\})\right)\mbox{]} \label{gamma b},
\end{eqnarray}
depending on which version of gauge $\g$ is in question.
 From this it is easy to find a cell or cells $I^r$ which form a $\g$-fine division of the cell $J^r$ in step 4 above, giving a contradiction. 
\item
Therefore, if 1 is valid, then for each point $z_T = (y_{t_1}, z_{t_2}, z_{t_3}, \ldots)$ of the hyperplane $Z^1_T$ the set $L\left(z_T\right)$ of the gauge $\g$ contains a co-ordinate label $t_k$ with $k>1$.
\item
Likewise, with step 5 completed, a ``hyper-hyper-plane'' $Z^2_T$ of the hyperplane $Z^1_T$ is considered, and the assumption of non-divisibility in 1 above implies that the set
\[
L(z^2_T) = L\left(y_{t_1}, y_{t_2}, z_{t_3},z_{t_4},\ldots \right)
\]
contains a co-ordinate label $t_k$ with $k>2$.
\item
But in step 6, a fixed point $y_T=(y_{t_1},y_{t_2},y_{t_3},\ldots)$ is arrived at. The set $L(y_T)$ is defined in the gauge $\g$, and consists of a finite set 
\[
\left\{s_1, s_2, \ldots ,s_n\right\} \subset T,
\]
with $s_1<s_2< \cdots < s_n =t_m$ for some $m$. So $k>m$ implies $t_k \notin L(y_T)$.
\item
Assumption 1 implies that iteration $m$ of the single-coordinate bisection process delivers a coordinate label $t_k$, with $k>m$, for which 
\[
t_k \in L(y_T)=L(y_{t_1},\ldots ,y_{t_m},y_{t_{m+1}}, y_{t_{m+2}},\ldots),\]

with $(y_{t_1},\ldots ,y_{t_m},y_{t_{m+1}}, y_{t_{m+2}},\ldots) \in Z^m_T$; contradicting step 11.

\end{enumerate}
Thus assumption 1 gives a contradiction, and must therefore be false.
The argument can be adapted for uncountable $T$, as demonstrated in the proof of Theorem 4 in \cite{MTRV}.
\nproof
\end{example}
The steps above establish divisibility for all three versions of gauge $\g$ in $\R^T$, including the version $\g^b$ defined in (\ref{definition of gauge g 0 (b)}). The condition labelled
(\ref{gamma b})
above is the critical step.
Thus the following corresponds to Theorem 4 (page 120--124 of \cite{MTRV}) and, as shown above, has almost identical proof.

\begin{theorem}
\label{gamma b divisibility of RT}
Given a gauge $\g$ (or $\g^b$) defined as (\ref{definition of gauge g 0 (b)}), there exists a $\g$-fine division of $\R^T$.
\end{theorem}
\proof
With minor adaptations, the proof follows that of Example \ref{divisibility infinite dim}.  \nproof

\section{Higher-dimensional Brownian Motion}
The reason for the preceding comments on product gauges is to seek to formulate three-dimensional Brownian motion in terms of the system of -complete integration. Returning to this issue, and following the scheme established in Example \ref{Bentropia}, let $T$ denote $]0, \tau]$, and let sample space $\Omega$ be
\[
R_{3T}=\prod \left\{\R^3: t \in T\right\}=\left(\R^3\right)^T.
\]
If $N$ is any finite subset of $T$, let $R_N$ denote 
$\prod \left\{\R^3: t \in N\right\}$, with finite-dimensional points $\bx_N$, cells $\bI_N$, association $(\bx_N,\bI_N)$ and gauges $\delta$ defined as in Example \ref{Bentropia} and (\ref{composite gauge for Rn}).
For any finite set $N \in \N(T)$, a cell $\bI=\bI_T =\bI_T[N]$ in $R_{3T}$ is
\[
\bI_T = \bI_N \times R_{T \setminus N}.
\]  
 So if $N = \{t_1, \ldots ,t_n\} \subset \,]0,\tau]$,
\[
\bI_T= \bI_T[N] = \left(\prod_{j=1}^n \left(I_1(t_j)\times I_2(t_j)\times I_3(t_j)\right) \right) \times R_{T \setminus N},
\]
where each component term of the Cartesian product $R_{T \setminus N}$ is $\R \times \R \times \R$. In physical terms, to say that $\bx_T \in \bI$ means that, for dimensions $i=1,2,3$ at time $t_j$ ($1 \leq j \leq n$), the $i$th co-ordinate of the spatial position of the particle lies in the real interval $I_i(t_j)$; and for any time $t \notin N$, the $i$th co-ordinate of the particle in space is unspecified, or arbitrary, for each of $i=1,2,3$.

A point-cell pair $(\bx_T, \bI_T[N])$ are \textit{associated} in $R_{3T}$ if $(\bx_N, \bI_N)$ are associated in $R_N$. That means that, as one-dimensional objects, $(x_i(t_j), I_i(t_j))$, $=(x_{i,j},I_{i,j})$, are associated in $\R$ for $i=1,2,3$ and for $t_j \in N$.

Next, define a gauge in $R_{3T}$, as follows.\footnote{Remember, $R_{3T}$ is $\left(\R^3\right)^{T}$, not $\R^{T}$.} For each $\bx_T, = \left(x_{1t},x_{2t},x_{3t}\right)_{t \in T},$ in $\bar R_{3T}$ (that is, $R_{3T}$ with points at infinity added), and for each finite set $N \in \N(T)$, let  
\[
L:  \left(\bar\R\times \bar\R\times \bar\R\right)^T  \rightarrow  \N(T), \mbox{ where }\; \bx_T \mapsto  L(\bx_T) \in \N(T), \;\mbox{ and }
\]
\begin{equation}\label{definition of delta N}
\delta_N: \left(\bar \R\times \bar \R \times \bar \R\right)^N
\rightarrow \left( \R_+\times \R_+ \times \R_+\right)^N 
\end{equation}
\[\mbox{ where }\;\bx_N  \mapsto \left(\delta_t(x_{1t}),\delta_t(x_{2t}),\delta_t(x_{3t})\right)_{t \in N}.
\]
A \textit{gauge $\g$ in} $R_{3T}$ consists of  
\begin{equation}
\label{definition of gauge g R3T} 
(L, \delta_N),\;\;N \supseteq L(x).
\end{equation} 
 Given a gauge $\g$, an associated point-cell pair $(\bx_T, \bI_T[N])$ in $R_{3T}$ is $\g$-\textit{fine} if $N \supseteq L(\bx_T)$ and if $(\bx_N,\bI_N)$ is $\delta_N$-fine in $R_N$. (The latter means that, for each $t \in N$, the one-dimensional pair $(x_i(t),I_i(t))$, $=(x_{i,t},I_{i,t})$, is $\delta_t$-fine in $\R$ for $i=1,2,3$ and for each $t \in N$.)

Armed with understanding from Example \ref{divisibility infinite dim}, a proof of $\g$-divisibility of the structured domain $R_{3T}, =(\R \times \R \times \R)^T,$ can be addressed as follows.
\begin{theorem}\label{divisibility of R3T}
Given a gauge $\g=(L, \delta_N)$ of $R_{3T}$, there exists a $\g$-fine division $\D$, $=$ $\{\left(\bx_T, \bI_T[N]\right)\}$, of $R_{3T}$.
\end{theorem}

\proof
Assume (for contradiction) that there is no $\g$-division of $R_{3T}$. The first lines (page 121 of \cite{MTRV}) of the proof of Theorem 4 carry forward unchanged. But the bisection in dimension $t_1$, indicated by cells $I^r$ at the bottom of page 121, need to be amended as follows. In dimension $t_1$ of $R_{3T}$, for $i=1,2,3$ let each of $I_i^r(t_1)$ denote one of the one-dimensional cells in the last line of page 121, and,  with additional subscript labels $1,2,3$ relating to the three Cartesian components of the single dimension $t_1$, let 
\[
I^r(t_1) = I_1^r(t_1) \times I_2^r(t_1) \times I_3^r(t_1),\;\;\;\mbox{ (or }\;\;
I^r_{t_1} = I_{1,t_1}^r \times I_{2,t_1}^r \times I_{3,t_1}^r),
\]
so $I^r(t_1)$ is a cell of $R^3$. (For emphasis the label $t_1$ can be in-line rather than subscripted.) Line 2 of page 122 in \cite{MTRV} therefore gives a fixed \[y(t_1) = (y_1(t_1),y_2(t_1),y_3(t_1))\; \in \; \bar R_{t_1} = \bar  \R \times \bar  \R \times \bar \R.\] The novel (and crucial) part of the proof occurs in lines 3 to 13 of page 122. For $R_{3T}$, this part of the proof carries forward practically unchanged. Likewise the re-definition of the gauge $\g$ in lines 14 to 19. The rest of the proof consists of iterations of these steps, leading to  contradiction of the initial assumption. \nproof

This result is the first step in providing a Riemann-observable domain or sample space $(\R \times \R \times \R)^T$ for three-dimensional Brownian motion. The proof for $n$-dimensional Brownian motion, $\left(\R^n\right)^T$ is essentially the same. And, if such is required, it works for domain
\[
\Omega = \prod \left(\R^{n_t}: t \in T\right)
\]
where $T$ is infinite and each $n_t$ is a positive integer depending on $t$. 
(In Theorem \ref{divisibility of R3T}, which is aimed at three-dimensional Brownian motion,  $n_t =3$ for each $t \in T$.)

{Theorem} \ref{divisibility of R3T} deals with divisibility of domain $\Omega=\left(\R \times \R \times \R\right)^T$ with $T$ infinite. The following Example describes a different structure, also involving infinite $T$ and  similar to $\left(\R^3\right)^T$.
\begin{example}\label{RT3}
Let
\[
\Omega = \R^T \times \R^T \times \R^T = \left(\R^T\right)^3,
\]
to be denoted $R_{T3}$ for short.
Elements of $R_{T3}$ have the form $\bx_T=\bx=(\bx_{1}, \bx_{2}, \bx_{3})$ where
\begin{equation}\label{RT3 el}
\bx_{} = \left(x_{T,1},x_{T,2}, x_{T,3}\right) =\left(x_{t,j}\right)_{t \in T,\;\;j=1,2,3}
=\left((x_{t,1}, x_{t,2}, x_{t,3})\right)_{t\in T},
\end{equation}
with $ x_{t,j} =x_{tj} \in \R$ for $j=1,2,3$. As to notation, it is sometimes easier to perceive the meaning if the label $j$ ($j=1,2,3$) is written as superscript instead of subscript. So the following representation is allowed:
\[
\bx=\bx_T =\left( \left(x_t^1 \right)_{t\in T},\left( x_t^2\right)_{t\in T}, \left( x_t^3\right)_{t\in T} \right).
\]
A cell $\bI_T$, $=\bI$, $=\bI[N]$ in $\left(\R^T\right)^3$ is 
$
I_1[N] \times I_2[N] \times I_3[N]
$
where $N = \{t_1, \ldots , t_n\}$ is any finite subset of $T$ and, for $t\in N$ and for $1 \leq j \leq 3$,
\[
I_j(t), = I_{tj}, = I_t^j,\; \in \I(\R);
\] 
so
\[
\bI = \bI_T[N]= 
\prod_{j=1}^3\left(\left( I_{t_1,j}\times I_{t_2,j} \times \cdots \times I_{t_n,j}\right) \times \R^{T \setminus N}\right).
\]
A pair $(\bx_T,\bI_T[N])$ are associated in $R_{T3}$ if, for each $t \in N$ and $j=1,2,3$, the one-dimensional pair $(x_{tj},I_{tj})$ are associated in $\R$.

\end{example}    
The new domain $\left(\R^T\right)^3$, $=R_{T3}$, is a (finite) product of (infinite) products.

The next step is to construct a gauge for the domain $R_{T3}$. The underlying method of construction of gauges for product spaces is to use already-defined gauges for the component domains of the product domain, and form a product of such gauges. This works even when the new gauge is a ``product of product-gauges''. So it is assumed that, for each $t \in T$ and each $y \in \bR$, one-dimensional gauges 
$\da_t(y)>0$ are given.

In the case of $\left(\R^T\right)^3$, for each element $\bx \in \left(\bar\R^T\right)^3$ let $L(\bx)$ be a finite subset of $T$; and, for each $N =\{t_1, \ldots , t_n\} \in \N(T)$ and for each \[
y = y_N = \left(\left(y_{t_i}^j\right)_{t_i \in N,\,j=1,2,3}\right)\;
\in \;\bar\R^N,
\]
 let 
\[
\da_N(y) =\left\{\left(\da(y_{t_i}^1) ,\da(y_{t_i}^2),\da(y_{t_i}^3)\right)_{t_i \in N}\right\}.
\]
Now define a gauge $\g$ in $\left(\R^T\right)^3$ as
\begin{equation}\label{definition of gauge g RT3} 
\g = \left(L, \da_N\right),
\;\mbox{
for all }\bx \in \left(\bar\R^T\right)^3\;\mbox{ and all }N \in \N(T).
\end{equation} 
Then the associated pair $(\bx_T, \bI_T[N])$ are $\g$-fine if $N\supseteq L(\bx_T)$ and $(x_{t,j},I_{t,j})$ are $\da_t$-fine in $\R$ for each $t \in N$. 

If $N$ contains $n$ elements, then a cell $\bI[N]$ has a total of $3n$ restricted one-dimensional component cells $I_{t,j}$ (or $I_t^j$) for $t \in N$ and for $j=1,2,3$. Essentially, a gauge for  $\left(\R^T\right)^3$ regulates each of these edges by means of the corresponding $\da$-values of the $3n$ associated points $y_t^j \in \bar \R$. The condition $N \supseteq L(\bx_T)$ in the gauge $\g$ ensures that the dimension sets $N$ that appear in Riemann sums contain arbitrarily large numbers of elements.

Here is a quick reminder of how the -complete integral of a function $h$ is defined by means of gauges. If $h(\bx_T, N, \bI[N])$ is a real- or complex-valued function, then $h$ is integrable on $\left(\R^T\right)^3$, with integral $\alpha = \int_{\left(\R^T\right)^3}h$, if the following holds. Given $\ve>0$, there exists a gauge $\g = (L, \da_N)$ such that, for every $\g$-fine division $\D$ of $\left(\R^T\right)^3$,
\[ 
\left|\alpha - (\D) \sum h(\bx_T, N, \bI[N]) \right| < \ve.
\]
This definition is empty and meaningless unless the required divisions $\D$ of $\left(\R^T\right)^3$ exists. When $\g$-fine divisions are shown to exist, elaboration of the properties of the -complete integral $\int h$ follows a common pattern, and the theory described in Chapter 4 of \cite{MTRV} is generally applicable.

Thus the basis of the integration in $\left(\R^T\right)^3$ is the following theorem.

\begin{theorem}\label{existence of g-division for RT3}
Given a gauge $\g$ for $\left(\R^T\right)^3$, as defined in (\ref{definition of gauge g RT3}), there exists a $\g$-fine division of $\left(\R^T\right)^3$.
\end{theorem}
\proof
This theorem is also valid if $\left(\R^T\right)^3$ is replaced by any cell $\bI$ of $\left(\R^T\right)^3$. Assume (for contradiction) that no $\g$-fine division of $\left(\R^T\right)^3$ exists. For simplicity, assume also that $T$ is a countable set,
$
T=\left\{\tau_1, \tau_2, \ldots \right\}
$,
so \begin{eqnarray*}
\left(\R^T\right)^3 &= &\R^T\times \R^T\times \R^T\vt
&=&\left(\R \times \R \times \R \times \cdots \right)\times
\left(\R \times \R \times \R \times \cdots \right)\times
\left(\R \times \R \times \R \times \cdots \right)\vt
&=&\left(\R \times \R^{T\setminus \{ \tau_1\}} \right)\times
\left(\R \times \R^{T\setminus\{ \tau_1\}} \right)\times
\left(\R \times \R^{T\setminus \{\tau_1\}} \right)
.
\end{eqnarray*}
(The proof can be adapted for uncountable $T$, as in Theorem 4 of \cite{MTRV}.)
The product domain has three components, corresponding to $j=1,2,3$; and each of the three components is a composite of factors
\[
\left(\R \times \R^{T\setminus \{ \tau_1\}} \right),\;=\; \left(\R^j_{\tau_1} \times \R^{T\setminus \{ \tau_1\}} \right)
\]
for $j=1,2,3$.
Now bisect, successively, the term $\R^j_{\tau_1}$ (jointly for $j=1,2,3$) of the three factors; so that, at each bisection, a non-$\g$-divisible cell is obtained. As in Example \ref{divisibility infinite dim},
the successive bisections yield 
\[
\left(\{y^1_{\tau_1}\} \times \R^{T\setminus \{ \tau_1\}} \right) \times 
\left(\{y^2_{\tau_1}\} \times \R^{T\setminus \{ \tau_1\}} \right) \times
\left(\{y^3_{\tau_1}\} \times \R^{T\setminus \{ \tau_1\}} \right).
\]
This procedure can be repeated successively for $j=2,3, \ldots$, leading to a succession of non-divisible cells in domains
\[
\left(\R^{T\setminus \{ \tau_1, \tau_2, \ldots , \tau_k\}} \right)
\]
for $k=2,3, \ldots$. Thus a point 
\[
y=\left(\left(y^1_{\tau_1},y^2_{\tau_1},y^3_{\tau_1}\right), \left(y^1_{\tau_2},y^2_{\tau_2},y^3_{\tau_2}\right), \left(y^1_{\tau_3},y^2_{\tau_3},y^3_{\tau_3}\right) \ldots \right) \in \left(\R^T\right)^3
\]
is arrived at by iteration. And, as in Example \ref{divisibility infinite dim}, if $L(y) = M$ and if $m = \max \{i: \tau_i \in M\}$, 
the original assumption of non-$\g$-divisibility fails at the $m$-th stage in the preceding iteration. \nproof

With corresponding division structures the preceding results can be established for domains $\left(\R^T\right)^n$ for any positive integer $n$. Also, variants of
(\ref{definition of gauge g RT3})
can be substituted, without causing difficulty in the preceding proof of $\g$-divisibility.

The domain $\left(\R^3\right)^T$ (or $R_{3T}$, with $T$ of infinite cardinality) was introduced in Example \ref{3-dim BM} as a sample space domain for three-dimensional Brownian motion. This has sample values $(x_{1t},x_{2t},x_{3t})$ for $t\in T$, representing the coordinates in $\R\times \R\times\R$ of a Brownian particle at any time $t$, $t \in T=\,]0,\tau]$.

In contrast, Example \ref{RT3} has domain $\left(\R^T\right)^3$, with sample values (or representative elements
\[
\left(x_{t1},x_{t2},x_{t3}\right)_{t\in T},
\]
as described in (\ref{RT3 el}).
Since the numerical value of $t$ is the same in each instance of $x_{jt}$ ($j=1,2,3$), there is no essential physical difference between Examples \ref{3-dim BM} and \ref{RT3}.

\section{Products of Product Spaces}
The progression has been from finite Cartesian product spaces $\R^n$ to infinite Cartesian product spaces $\R^T$; and then to  compound finite-infinite Cartesian product domains. The latter were described in Examples \ref{3-dim BM} and \ref{RT3}; respectively,
\[
\left(\R^3\right)^T,\;\;\;\left(\R^T\right)^3; \;\;\;\;\mbox{ or }\;\;\left(\R^P\right)^T,\;\;\;\left(\R^T\right)^P,
\]
with $P$ finite and $T$ infinite.  It is natural to inquire whether the -complete structure can be formulated for $\Omega =\left(\R^P\right)^T$ with both $T$ and $P$  infinite (for simplicity, countably infinite).  
Before tackling this, the following example gives a different interpretation of ${\R^T}^3$.
\begin{example}\label{RT.3}
Suppose $T^j$ are (countably) infinite labelling sets for $j=1,2,3$, and suppose $\Omega$ is a composite infinite-finite Cartesian product domain,
\begin{equation}\label{RT1T2T3}
\Omega = \R^{T_1} \times \R^{T_2} \times \R^{T_3},
\end{equation}
whose elements are 
\begin{equation}\label{RT1T2T3 el}
\bx =\bx_\Omega= \left(x_{T_1}, x_{T_2}, x_{T_3}\right)
=\left((x_{t^1})_{t^1 \in T_1}, (x_{t^2})_{t^2 \in T_2}, (x_{t^3})_{t^3 \in T_3}\right).
\end{equation}
If $T_1=T_2=T_3 =T$, then $\Omega$ has the form $\R^T \times \R^T \times \R^T$; but not with the same meaning as the space $R_{T3}$ above. If we follow (\ref{RT1T2T3 el}), a representative element of $\Omega$ is  
\begin{equation}\label{RTTT el}
\bx =\bx_\Omega= \left(x_{T,1}, x_{T,2}, x_{T,3}\right)
=\left((x_{t,1})_{t \in T}, (x_{t,2})_{t \in T}, (x_{t,3})_{t \in T}\right).
\end{equation}
Despite superficial appearances, the space of the latter elements is \textbf{not} the same as the space $R_{T3}$ of Example \ref{RT3} which uses the symbolic name $\left(\R^T\right)^3$. The difference between these two domains is expressed in the difference between their representative elements, in (\ref{RT3 el}) and (\ref{RTTT el}). In (\ref{RT3 el}) the $t$ of $x_{t,j}$ is the same for $j=1,2,3$, while in (\ref{RTTT el}) the $t$ of $x_{t,j}$ may be different for $j=1,2,3$. To make the distinction clearer, denote the space consisting of elements (\ref{RT3 el}) and (\ref{RTTT el}) by, respectively, 
\[
\left(\R^{T}\right)^3,\;\;\;\;  \R^{\left(T^3\right)},
\]
or, for short, $R_{T3}$,   $R_{T.3}$, respectively. Each element (\ref{RT3 el}) is also an element (\ref{RTTT el}), but not conversely; so $\left(\R^{T}\right)^3\subset \R^{\left(T^3\right)}$.

Cells $\bI$ in $R_{T.3}$ are constructed as Cartesian products of cells 
\[
I^j[N^j] \in \I(\R^T),\;\;\; j=1,2,3,
\]
as follows. With $N^j \in \N(T)$, $j=1,2,3$, let
\[
\bN=\left(N^1,N^2,N^3\right).
\]
Write this as $\bN \in \N(T \times T \times T)$, and define
\[
\bI=\bI_{T \times T \times T}=\bI_{T.3}=\bI[\bN] := I^1[N^1] \times I^2[N^2]\times I^3[N^3];
\]
so
\begin{equation}\label{RT.3 cell}
\bI= \bI_{T\times T \times T} =\bI_{T.3} = \bI[\bN] =
\prod_{j=1}^3 \left( \prod_{t \in N^j} I_{t,j} \times \R^{T \setminus N^j} \right),
\end{equation}
where, for $j=1,2,3$ and for $t \in N^j$, $I_{t,j}$ (or $I_t^j$) is\footnote{As to notation, the label $j$ is used in both subscript and superscript form to distinguish between the component domains $\R^T$ in $\R^T \times \R^T \times \R^T$ (or $\R^{\left(T^3\right)}$, or $R_{T.3}$). Accordingly,   it is plausible to write the one-dimensional component cells   $I_{t,j}$ as $I_t^j$, depending on which version provides  better intuition.} a one-dimensional real interval, $I_{tj}=I^j_t \in \I(\R)$.

A point-cell pair $(\bx,\bI[\bN])$ are associated in $R_{T.3}$ if, for $j=1,2,3$, 
\[
(x_{T,j},I^j[N^j])\mbox{ are associated in }\R^{T_j}, =\R^T;
\] so $(x_{t_j}^j,I_{t_j}^j)$ are associated in $\R$ for each $t_j\in N^j$, $j=1,2,3$.

Next, define a gauge $\g$ in $R_{T.3}$. Assume that, for $j=1,2,3$ and each $t \in T$ there is a gauge $\da_t^j(y)>0$ for $y \in \bar \R$. 
Let 
\[
\bL(\bx)=\left(L^1(x^1_T), L^2(x^2_T),L^3(x^3_T)\right) \in \N(T)\times \N(T)\times \N(T).
\]
A gauge $\g$ in $R_{T.3}$ is 
\begin{equation}
\label{definition of gauge g RT.3} 
\g=\left(\bL,\da\right)\;\mbox{ where }\; \da = \left\{\da_t^j\right\}_{t \in T,\; j=1,2,3}.
\end{equation}
With $\bN=(N^1,N^2,N^3) \in \N(T)\times \N(T)\times \N(T)$, an associated pair $(\bx,\bI[\bN])$ is $\g$-fine in $R_{T.3}$ if, for $j=1,2,3$, \[ N^j \supseteq L^j(x^j_T)\mbox{ and
}(x_{t_j}^j,I_{t_j}^j)\mbox{ is }\da_{t_j}^j\mbox{-fine in }\R
\] for each $t_j\in N^j$. Existence of $\g$-fine divisions of $R_{T.3}$ is proved in Theorem \ref{existence of g-fine divisions of RT.3} below.

\end{example} 
\begin{theorem}\label{existence of g-fine divisions of RT.3}
Given a gauge $\g=(\bL,\da)$, there exists a $\g$-fine division of $R_{T.3}$.
\end{theorem}
\proof
For each $j$ ($j=1,2,3$), take $\g^j = (L^j, \da^j)$ where $\da^j= \{\da_t\}_{t \in T}$. Then, by the construction in (\ref{definition of gauge g RT.3}), there exists a $\g$-fine division of $R_{T.3}$ if and only there exists a $\g^j$-fine division of $\R^T$ for each $j$. But the latter holds by Theorem 4 of \cite{MTRV}.
\nproof

\section{Illustration of Products of Products}
In order to distinguish clearly between the two domains $R_{T3}$ and $R_{T.3}$, 
and to provide motivation, clarification and intuition,
here are some illustrations which demonstrate the different roles of $R_{T3}$ and $R_{T.3}$. 

\begin{enumerate}
\item[$\{R_{T3}\}$]\label{RT3 illustration}
Suppose $T=\;]0,\tau]$, and, for $t \in T$, $x_{tj}$ (or $x_t^j$) is the share value of one of three stocks, labelled $j$ for $j=1,2,3$. If, for each of the three stocks, and for any given time $t$, the share value $x_{tj}$ has a distribution function $F_{X_{tj}}$ defined on cells $I_{tj} \in \I(\R)$ then, for each $t,j$, the (unpredictable) share value $x_{tj}$ is a sample value of an observable 
\begin{equation}\label{RT3 rv basic 1}
X_{tj} \sq x_{tj}\left[F_{X_{tj}}, \R\right].
\end{equation}
The general real-valued ample space $\Omega = \R$ is assumed.) If, for each $j$, there is a joint distribution function $F_{X_{Tj}}$ defined on cylindrical intervals $I[N] \in \I(\R^T)$ then $X_{Tj}$ (or $X_T^j$) is a joint observable (or joint process)
\begin{equation}\label{RT3 rv joint basic 1}
X_{Tj} \sq x_{Tj}\left[F_{X_{Tj}}, \R^T\right].
\end{equation}
Now, for any given $t \in T$, consider the joint share values 
\begin{equation}\label{RT3 joint share value}
\bx_t = \left(x_{t1}, x_{t2}, x_{t3}\right), \,= \left(x_{t}^1, x_{t}^2, x_{t}^3\right).
\end{equation} 
Suppose, for each $t$ and for $j=1,2,3$, there exists a joint distribution function $F_{\bX_t}$ defined on cells 
\begin{equation}\label{RT3 cell}
\bI_t = I_{t1}\times I_{t2}\times I_{t3}, = I_{t}^1\times I_{t}^2\times I_{t}^3 \in \I\left(\R^3\right).
\end{equation}
Suppose, for each finite set $N =\{t_1, \ldots , t_n\}\in \N(T)$ and 
\begin{equation}\label{RT3 joint}
\bI_T[N] =  \prod_{j=1,2,3}\left(I_{t_1 j} \times  I_{t_2 j} \times \cdots  \times I_{t_n j}  \times \R^{T\setminus N}\right),
\end{equation}
there is a joint distribution function 
$F_{\bX_T}(\bI_T[N]$ defined on cells $\bI_T[N] \in \I\left(R_{T3}\right)$, then $\bX_T$ is an observable process
\begin{equation}\label{RT3 joint observable}
\bX_T \sq \bx_T \left[F_{\bX_T} , \left(\R^T\right)^3 \right].
\end{equation}
Now suppose a derivative asset $Y, =Y_\tau$, is constructed so that at time $\tau$ the holder of the derivative asset receives a payment of amount $\kappa$ if, at any time $t \in T=\,]0,\tau]$, the sum of the values of the three underlying assets exceeds a specified  amount $\lambda$; and otherwise the holder of the derivative receives nothing. Thus
$Y = f(\bX_T)\;=$
\begin{equation}\label{RT3 example}
= f(X_{T1},X_{T2},X_{T3}) = \left\{
\begin{array}{ll}
\kappa & \mbox{if } \; 
\sum_{j=1}^3 x_{t,j}
> \lambda \mbox{ for some } t \in T,\vt
0 & \mbox{otherwise.}
\end{array}
\right.
\end{equation}
Then the contingent observable $Y$ is a contingent random variable if its expected value  exists:
\begin{equation}\label{RT3 expectation}
\E[Y] = \int_{R_{T3}} f(\bX_T) F_{\bX_T} \left(\bI[N]\right).
\end{equation}
The latter (if it exists) is the expected payoff of the derivative asset $Y$.
\item[$\{R_{T.3}\}$]\label{RT.3 illustration}
The next step is to provide a practical illustration of joint random variation on domain
\[
\Omega = \R^{\left(T^3\right)} = R_{T.3},  
\]
and, in doing so, to take note of how this kind of random variation contrasts  with the preceding. As before,  $x_{tj}$ (or $x_{t}^j$) is the share value of any one of three stocks,  ($j=1,2,3$). As before, a distribution function $F_{X_{tj}}$ defined on cells $I_{tj} \in \I(\R)$ is assumed, with  observable 
\begin{equation}\label{RT.3 rv basic 1}
X_{tj} \sq x_{tj}\left[F_{X_{tj}}, \R\right].
\end{equation}
Also, a joint distribution function $F_{X_{Tj}}$ is defined on cylindrical intervals $I[N] \in \I(\R^T)$, so $X_{Tj}$ is a joint observable (or joint process)
\begin{equation}\label{RT.3 rv joint basic 1}
X_{Tj} \sq x_{Tj}\left[F_{X_{Tj}}, \R^T\right].
\end{equation}
So far, everything is as before. But at this point we depart from the previous scheme. In the preceding example $\bx_t$ represented  $\left(\bx_{t,j}\right)_{j=1,2,3}$, or $\left(x_{t,1}, x_{t,2}, x_{t,3}\right)$, in which time $t$ is the same for each component.
But what is now needed is joint share values of the form described in Example \ref{RT.3},    
\begin{equation}\label{RT.3 joint share value}
\bx_{T \times T \times T}, = \bx_{T^3}, = \left(x_{t_1,1}, x_{t_2,2}, x_{t_3,3}\right)_{t_1,t_2,t_3 \in T}
, = \left(x_{t_1}^1, x_{t_2}^2, x_{t_3}^3\right)_{t_1,t_2,t_3 \in T},
\end{equation} 
where the times $t_j$ may be different for $j=1,2,3$. Thus 
\[
\bx_{T^3} = \left(x_{T,1},x_{T,2},x_{T,3}\right);\;\;\;\mbox{ or }\;\bx_{T^3} =\left(x_T^1,x_T^2,x_T^3\right),
\]
where the latter representation emphasises the fact that $R_{T.3}$ is the Cartesian product $\R^T \times \R^T \times \R^T$. Write $\bx_{T^3}$ as $\bx$ (or as $\bx_{T.3}$ if clarity requires this). 
As described in Example \ref{RT.3}, cells $\bI$ in $R_{T.3}$ are  Cartesian products of cells 
\[
I^j[N^j] \in \I(\R^T),\;\;\; j=1,2,3.
\]
As in Example \ref{RT.3}, with $N^j \in \N(T)$ and $j=1,2,3$, 
\begin{eqnarray}
\bN&=&\left(N^1,N^2,N^3\right)\in \left(\N(T)\times \N(T) \times \N(T)\right),\nonumber \vt
\bI&=&\bI_{T \times T \times T}=\bI_{T.3}=\bI[\bN] := I^1[N^1] \times I^2[N^2]\times I^3[N^3], \nonumber\vt
\label{RT.3 cell+}
&=& 
\prod_{j=1}^3 \left( \prod_{t \in N^j} I_{t,j} \times \R^{T \setminus N^j} \right),
\end{eqnarray}
where, for $j=1,2,3$ and for $t \in N^j$, $I_{t,j}$ (or $I_t^j$) is a one-dimensional real interval, $I_{tj}=I^j_t \in \I(\R)$.


\end{enumerate}

\textit{\textbf{Now suppose there is some joint measurement or   observation of unpredictable values $x_{t,j}$ (or $x_t^j$) for $t \in T$ and $j=1,2,3$, and suppose these unpredictable values are regulated by
 a joint distribution function 
\begin{equation}\label{RT.3 joint dist fn}
F_{\bX_{T.3}}(\bI[\bN]\;\mbox{ defined on cells }\;\bI[\bN] \in \I\left(R_{T.3}\right).
\end{equation} 
Then $\bX_{T.3}$ is an observable process
\begin{equation}\label{RT.3 joint observable}
\bX_{T.3} \sq \bx_{T.3} \left[F_{\bX_T} , \R^{\left(T^3\right)} \right].
\end{equation}}}
The next step is to envision some practical scenario which realizes this abstract conceptual framework of domains (or sample spaces) and observables. An adaptation of the example given in $\{RT3\}$ above provides illustration.

Accordingly, suppose a derivative asset $Y, =Y_\tau$, is constructed so that at time $\tau$ the holder of the derivative asset receives a payment of amount $\kappa$ if, at any time $t^j$ ($0<t^j\leq \tau$, $j=1,2,3$), the sum of the values of the three underlying assets exceeds a specified  amount $\lambda$; and otherwise the holder of the derivative receives nothing. Thus
\begin{equation}\label{RT.3 example}
Y = f(\bX_{T.3})
= \left\{
\begin{array}{ll}
\kappa & \mbox{if } \; 
\sum_{j=1}^3 x^j_{t^j}
> \lambda \mbox{ for some } t^1,t^2,t^3 \in T,\vt
0 & \mbox{otherwise.}
\end{array}
\right.
\end{equation}
Then the contingent observable $Y$ is a contingent random variable if its expected value  exists:
\begin{equation}\label{RT.3 expectation}
\E[Y] = \int_{R_{T.3}} f(\bX_{T.3}) F_{\bX_{T.3}} \left(\bI[\bN]\right).
\end{equation}
The latter (if it exists as an integral value) is the expected payoff of the derivative asset $Y$.

\textit{\textbf{The difference between (\ref{RT.3 example}) and (\ref{RT3 example}) is that, in the latter the time $t$ is the same for the three underlying assets $j=1,2,3$, while in the former the three times $t^1,t^2,t^3$ may be different for some/all of the three underlying assets.}} 

\section{Infinite Products of Infinite Products}
\begin{example}\label{RP.T}
It is straightforward to extend the meaning of Examples \ref{RT3} and \ref{RT.3} to domains
\[
\left(\R^T\right)^n = R_{Tn} \;\;\mbox{ and }\;\;\R^{\left(T^n\right)} = R_{T.n},
\]
where $n$ is any positive integer.

In the development of the theory, finite-dimensional domains $\R^n$ have been extended to infinite-dimensional domains $\R^T$.
This suggests that, in $\left(\R^T\right)^n$ and/or $\R^{\left(T^n\right)}$, the finite  $n$ may be replaced by an infinite set $P$ of dimension-labels, giving
\[
\left(\R^T\right)^P, = R_{TP} \;\;\mbox{ and/or }\;\;\R^{\left(T^P\right)}, =R_{T.P}.
\] 
We will now provide some meaning to these; specifically the latter, $R_{T.P}$. For simplicity it is assumed that both $T$ and $P$ are countable.

Since $P = \{p_1,p_2, \ldots \}$ is countable the members of the set $\Omega = \R^{\left(T^P\right)}, =R_{T.P},$ consist of sequences $\left(x_{p_1}, x_{p_2}, \ldots, x_{p_j}, \ldots \right)$, where each term 
$x_{p_j}$ is itself a sequence 
\[
\left(x_{t_1, p_j}, x_{t_2,p_j}, \ldots , x_{t_i,p_j}, x_{t_{i+1},p_j}, \ldots \right),\;\;\;\mbox{ or }\;\;
\left(x_{{1},{ j}}, x_{{2},{j}}, \ldots , x_{{i},{j}}, x_{{{i+1}},{j}}, \ldots \right),
\]
or 
$\left(x_{1}^{ j}, x_{2}^{j}, \ldots , x_{i}^{j}, x_{{i+1}}^{j}, \ldots \right)$
where these three equivalent notations are used in various contexts. Thus, if $x \in \Omega= R_{T.P}$,
\[
\bx, = \bx_{T.P},= \left(\left(x_{ij}\right)_{t_i \in T}\right)_{p_j \in P},
\]
and, in   notation whose ambiguity is exemplified in (\ref{RT3 example}), (\ref{RT.3 example}), the domain $\Omega$ is
\[
\R^{ \left(T^P\right)} =\prod_p \left(\prod_t\left(\R:t \in T\right): p \in P\right).
\]
The next step is to specify cells $\bI$ in $\Omega$. The domain $ \Omega$ can be represented as $S^P$ with $S=\R^T$. For $N_P = \{p_1, \ldots , p_n\} \in \N(P)$, a cell $\bI$ of $S^P$ is
\[
\bI =\left(J_{p_1} \times \cdots \times J_{p_n}\right) \times S^{P \setminus N_P}, = J(N_P) \times \left(\R^T\right)^{P \setminus N_P},
\]
where each $J_{p_j}$ is a cell $\bI_{T,j}$ (or $\bI_T^j$) in $\R^T$. So if 
\[
N_{T,j}=N_T^j= \{t_{1j}, \ldots , t_{m_j,j}\}, =\{t_{1}^{j}, \ldots , t_{m_j}^{j}\}, \in \N(P),
\] a cell $J_{p_j}$ of $\R^T$ is 
\begin{eqnarray*}
J_{p_j} &=& \bI_{T,j} \;=\; I(N_{T,j}) \times \R^{T \setminus N_{T,j}}\vt
&=& I_{t_{1},p_j} \times I_{t_{2},p_j} \times \cdots \times I_{t_{m_j},p_j} \times \R^{T \setminus N_{T,j}},\vt
& =& I_{t_{1}}^{p_j} \times I_{t_{2}}^{p_j} \times \cdots \times I_{t_{m_j}}^{p_j} \times \R^{T \setminus N_{T,j}}.
\end{eqnarray*}
Thus, for $N_P \in \N(P)$ and $N_{T,j} \in \N(T)$, a cell of domain $\Omega=R_{T.P}$ is
\begin{eqnarray}\label{cell RPT}
\bI_{T.P} &= & \left(\prod_{j=1}^n \bI_T[N_{T,j}]\right) \times {\R^T}^{P \setminus N_P} \vt
&=&\left(\prod_{j=1}^n\left(
I_{t_{1},p_j} \times I_{t_{2},p_j} \times \cdots \times I_{t_{m_j},p_j} \times \R^{T \setminus N_{T,j}}
\right)\right)
\times {\R^T}^{P \setminus N_P}
\nonumber\vt
& = &\left(\prod_{p_j \in N_P} \left( \left( \prod_{t_{ij} \in N_{T,j}} I_{t_{ij},p_j}\right) \times \left(\prod_{t \in T\setminus N_{T,j}}\R\right) \right)\right)
\times \left(\prod_{p \in P\setminus N_P}\R^T\right).
\nonumber
\end{eqnarray}
Loosely, $\bI_{T.P}$ can be written $\bI\left[N_{T_P}\right]\left[N_P\right]$.

Consider a point $\bx=\bx_{T.P}$ 
and a cell $\bI=\bI_{T.P}$  of ${\R^T}^P$, so
\begin{equation}\label{RT.P (x,I)}
 \bx_{T.P}= \left(\left(x_{t,p}\right)_{t \in T}\right)_{p \in P},\;\;\;\;
\bI_{T.P} = \left(\prod_{j=1}^n \bI_T[N_{T,j}]\right) \times {\R^T}^{P \setminus N_P},
\end{equation}
with $N_P =\{p_1, \ldots ,p_j, \ldots , p_n\}$ and $N_{T,j} =\{t_{1j}, \ldots , t_{m_j,j}\}$ for $j=1, \ldots ,n$. Then
the pair $(\bx_{T.P},\bI_{T.P})$ are associated in $R_{T.P}$ if, for each $N_{T,j}$ ($j=1, \ldots , n$), the cell $\bI(N_{T,j})$ has the point $\bx(N_{T,j})$ as a vertex in $\R^{N_{T,j}}$. 
Expressed as points and intervals of $\R$, this means that $x_{t_i,p_j}=x_{t_i}^{p_j} \in \bar\R$ is a vertex of the corresponding restricted cell $I_{t_i,p_j} = I_{t_i}^{p_j}\subset \R$ for ${1} \leq {i} \leq {m_j}$, $1 \leq j \leq n$.
\end{example}

To define -complete integrals on ${\R^T}^P$, gauges must be formulated to regulate the associated pairs $(\bx_{T.P}, \bI_{T.P}[N_{T_P}][N_P])$ used in forming Riemann sums. It is assumed that, for each component $\R$ of the compound product space ${\R^T}^P$ there is a gauge $\da_t^p(y)>0$ for $y \in \bar\R$; and it is assumed that for each component product space $\R^T$, ($=\R^{T,p}$, $p \in P$), of ${\R^T}^P$, there is a gauge $\g^p =(L^p, \{\da_t^p\})$ in $\R^T$, constructed in the standard way, for standard product space $\R^T$.

A gauge $\g $ for the compound product space ${\R^T}^P$ is formed from the component gauges $\g^p$ and $\da_t^p$, as follows. For each $\bx$ of 
${{\overline{\R}}^T}^P$, let $\bL(\bx)$ be an element $N=N_P \in \N(P)$, and let
\[
\g=\g_{T.P} := \left(\bL, \left\{\g^p\right\}_{p \in P}, \left\{\da_t^p\right\}_{t\in T,\;p\in P}\right).
\]
Then  the associated pair $(\bx, \bI) = (\bx_{T.P}, \bI_{T.P})$ of $R_{T.P}$ is $\g$-fine if, in (\ref{RT.P (x,I)}),
\begin{eqnarray}\label{RT.P g-fine}
N_P&=&\{p_1, \ldots , p_n\} \supseteq \bL(\bx), \nonumber\vt
\left(x_{T}^{p_j}, I_{T}^{p_j}\right) &\mbox{is}& \g^{p_j}\mbox{-fine for } j=1, \ldots n. 
\end{eqnarray}
In other words, $(\bx,\bI)$ is $\g$-fine if
\begin{eqnarray}
\{p_1, \ldots , p_n\} & \supseteq & \bL(\bx), \label{RT.P g-fine (1)}\vt
\{t_1^{p_j}, \ldots , t_{m_j}^{p_j}\} &\supset  &L^{p_j}(x_T^{p_j}) \mbox{ for } j=1, \ldots ,n, \;\;\;\mbox{ and}\label{RT.P g-fine (2)}\vt 
(x_{t_i}^{p_j}, I_{t_i}^{p_j}) & \mbox{is}& \da_{t_i}^{p_j}\mbox{-fine in } \R \mbox{ for } 1\leq i \leq m_j,\;\;1 \leq j \leq n. \label{RT.P g-fine (3)}
\end{eqnarray}
The purpose of the gauge construction is to define the -complete integral of some function $h(\bx,\bI)$ on the domain $\R^{T^P}$. For this to work there have to be $\g$-fine divisions of the domain in order to be able to form Riemann sum estimates of the integral.
\begin{theorem}\label{RT.P existence of g-fine divisions}
Given a gauge $\g$ in $\R^{T^P}$, there exists a $\g$-fine division of $\R^{T^P}$.
\end{theorem}
\proof
The general idea of the proof is a development or extension of the method of successive bisection combined with successive slicing described in Example \ref{divisibility infinite dim}. Assume, for contradiction, that there is no $\g$-fine division of the domain. For simplicity the sets $T$ and $P$ are assumed countable, so write
$
P=\{p_1,p_2, \ldots,\}$, $T=\{t_1,t_2, \ldots \}$.
Select dimension $p_1$ so, with $S=\R^T$ and $P_1 = P\setminus \{p_1\}$,
\[
{\R^T}^P = S \times\left( \left({\R^T}\right)^{P\setminus \{p_1\}}\right), =S\times R_{T.P_1}.
\]
Successively bisect $S=\R^{\{t_1,t_2, \ldots \}}$ as follows. For $s\;=$
\[
-q2^q,\;\;\;\;-q2^q+1,\;\;\;\;-q2^q+2, \;\;\;\ldots \;,\;\;\;\; q2^{-q}-1,
\]
write $I_{sq}$ as one of the $(s,q)$-binary one-dimensional cells
\[
]-\infty,-q],\;\;\; \;\;\;\;\;\;]s2^{-q},\;(s+1)2^{-q}],\;\;\;\;\;\;\;\;\;]q,\;\infty[.
\]
For $r=1,2,3, \ldots$, let $N_r=\{t_1, t_2, \ldots , r\}$, and, for each $q$ let $J^{qr}$ be a cell in $S=\R^T$,
\[
J^{qr} = J^{qr}[N_r] = \left(I_{1,q} \times I_{2,q} \times \cdots \times I_{rq} \right) \times \R^{T \setminus N_r}
\]
where, for $s=1, 2, \ldots , r$, $I_{sq}$ is an arbitrary $(s,q)$-binary one-dimensional cell. For each $r,q$, the cell $J^{qr}$ in $S=\R^T$ is divisible with respect to the product gauge formed from (\ref{RT.P g-fine (2)}) and (\ref{RT.P g-fine (3)}) above. (This follows from Theorem 4 of \cite{MTRV}.) For $r,q=1,2,3, \ldots$, the intersection of the closures of $J^{qr}$ contains a point $y_1 \in S=\R^T$. Consider the points
\[
\left(y_1, \bx_{T.P_1}\right)\;\;\mbox{ for }\;\;\bx_{T.P_1} \in R_{T.P_1}=\left({\R^T}\right)^{P\setminus \{p_1\}}.*
\]
If each of $\bL\left(y_1, \bx_{T.P_1}\right)$ consists only of $p_1$ then (as in Theorem 4 of \cite{MTRV}, or Example \ref{divisibility infinite dim} above), the assumption of non-$\g$-divisibility of the domain $R_{T.P}$ fails. Accordingly, as in Theorem 4 and Example \ref{divisibility infinite dim}, repeat the argument for $p_2,p_3, \ldots$ so that a point
\[
\mathbf{y}= \left(y_1,y_2,y_3, \ldots \right) \in R_{T.P} =\left({\R^T}\right)^{P}
\]
is obtained, for which $\bL(\mathbf{y})$ is an infinite set. This is impossible, so the original assumption of non-$\g$-divisibility of $R_{T.P}$ is disproved. \nproof

As in Theorem 4, this is extendible to uncountable sets $T$ and $P$. The proof can be adapted to domains 
\[
{{{\R^T}^P}^{\ldots}}^{Q} ,
\]
where the labelling sets $T,P, \ldots , Q$ are infinite.

\section{Integration in Product Spaces}
A gauge enables us to produce a system of integration in the domain.
The preceding examples illustrate  gauges in various product spaces.

But what about the corresponding integrals?
No examples, properties, theorems, proofs have been given above for any integrals.
Fortunately, the required integration theory, including variation theory, Fubini's theorem, and limit theorems, follows directly from the abstract development given in \cite{Henstock Lecture Notes}, and in \cite{H5,H6,H7,H4}.

The basic principles of this abstract \textit{division system}, or Henstock integral, are outlined in the Axioms of \cite{MTRV}, pages 111--113. The fundamental point is DS3 (the \textit{Division Axiom}), which has been established for the various product domains above. Because the cells $I$ of these product spaces are essentially rectangular, the other axioms (DS1, DS2, DS4--DS8) apply in a straightforward way.

Therefore, since the abstract proofs are already available, it is not necessary to re-hash them here for the various product space examples.

\[
\;
\]
{$\mbox{ }$ }\hfill Pat Muldowney, 1 December 2017, 
{$\mbox{ }$ } \hfill{\small\texttt{stieltjes.complete.integral@gmail.com}}


\begin{thebibliography}{99}



\bibitem{H5} Henstock, R., \emph{Generalised integrals of vector-valued functions},
Proceedings of the London Mathematical Society 19 (1969), 509--536.
\bibitem{H6} Henstock, R., \emph{Integration in product spaces, including Wiener and Feynman integration},
Proceedings of the London Mathematical Society 27 (1973), 317--344.
\bibitem{H7} Henstock, R., \emph{Integration, variation and differentiation in division spaces},
Proceedings of the Royal Irish Academy 78A(10) (1978), 69--85.
\bibitem{H4} Henstock, R., \emph{The General Theory of Integration}, Clarendon, Oxford, 1991.
\bibitem{Henstock Lecture Notes} Henstock, R.,  \emph{Lecture Notes} 
  1970, 70--132, \newline
  {\small\texttt{https://arxiv.org/pdf/1602.02993.pdf}}
\bibitem{MTRV} Muldowney, P., \emph{A Modern Theory of Random Variation}, John Wiley \& Sons Inc., Hoboken, New Jersey, 2012,
 \newline
  {\small\texttt{https://sites.google.com/site/stieltjescomplete/}}
  

  
\end{thebibliography}
\end{document}